\documentclass[12pt]{article}
\usepackage{amsmath}
\usepackage{amsfonts}

\advance \textheight by \topskip
\advance \textheight by 1pt
\makeatother
\def\bea{\begin{eqnarray}}
\def\ena{\end{eqnarray}}

\def\lar{\longrightarrow}
\def\non{\nonumber}

\def\deg{\hbox{deg}}

\def\dim{\hbox{dim}}

\def\ker{\hbox{Ker}}
\def\coker{\hbox{Coker}}
\def\ind{\hbox{index}}

\def\ord{\hbox{ord}}
\def\Res{\hbox{Res}}
\def\homega{\widehat{\omega}}

\def\tilp{\tilde{p}}

\def\tilp0{\tilde{p}_0}

\def\lar{\longrightarrow}

\def\qh{\widehat{q}}
\def\bxi{\mbox{\boldmath $\xi$}}

\newcommand{\bc}[2]{
\left(
\begin{array}{c}{#1}\\{#2}\end{array}
\right)}
\newcommand{\qbc}[2]{
\left[
\begin{array}{c}{#1}\\{#2}\end{array}
\right]}

\newcommand{\qed}{\hbox{\rule[-2pt]{3pt}{6pt}}}
\newtheorem{prop}{Proposition}
\newtheorem{theorem}{Theorem}
\newtheorem{defn}{Definition}
\newtheorem{lemma}{Lemma}
\newtheorem{cor}{Corollary}

\title{
Sigma Function as  A Tau Function
}

\author{
Atsushi Nakayashiki\thanks{
e-mail: 6vertex@math.kyushu-u.ac.jp}\\
Department of Mathematics, Kyushu University\\
}

\date{}
\begin{document}
\maketitle

\begin{abstract}
The tau function corresponding to the affine ring of a certain 
plane algebraic curve, called $(n,s)$-curve, embedded in the 
universal Grassmann manifold is studied. It is neatly expressed 
by the multivariate sigma function. This expression is in turn 
 used to prove fundamental properties on the series expansion 
of the sigma function established in a previous paper 
in a different method.

\end{abstract}

\section{Introduction}
The purpose of this paper is to study the multivariate sigma function
associated to certain plane algebraic curves, called (n,s)-curves, by means
of the tau function of the KP-hierarchy.

An (n,s)-curve is a plane algebraic curve given by the equation,
\bea
&&
y^n=x^s+\sum_{si+nj<ns}\lambda_{ij}x^iy^j,
\non
\ena
where $n,s$ are coprime and satisfy $1<n<s$. 
The sigma function associated to an (n,s)-curve had been introduced by
Buchstaber-Enolski-Leykin \cite{BEL1,BEL2} extending the Klein's 
hyperelliptic \cite{K1,K2} and Weierstrass' elliptic sigma functions. 
It is defined by modifying 
the Riemann's the theta function in such a way that it becomes 
modular invariant \cite{BEL1}.
If the genus of the curve is $g$, 
it is a holomorphic function of $g$ variables and has a remarkable
algebraic properties. Namely its series expansion at the origin begins at
the Schur function associated to the gap sequence at $\infty$ 
and all the coefficients of the expansion are homogeneous polynomials 
of $\{\lambda_{ij}\}$ with respect to certain degree. 
These properties have been proved in \cite{N}
by making an expression of the sigma function in terms of algebraic integrals
generalizing the Klein's formula \cite{K1,K2}. 
We remark that the general terms of the expansion are not known 
explicitly except for the elliptic case \cite{W} and the case of genus two 
\cite{BL}, where the recursion relations among expansion coefficients 
are explicitly given. Linear differential equations satisfied by 
sigma functions have been constructed in \cite{BL,BL08}. It can be 
a base to study the series expansions in more general cases. 
In this paper we propose another observation on the series expansions 
on sigma functions.

It is well known that the Weierstrass' elliptic functions 
$\wp(u)$ and $\wp'(u)$ uniformize the family of elliptic curves
$y^2=4x^3-g_2x-g_3$ even when $(g_2,g_3)$ is contained in the discriminant.
Similarly the second and third logarithmic derivatives of the Klein's
hyperelliptic sigma function uniformize the family of affine Jacobians, 
that is, the Jacobian minus the theta divisor \cite{BEL1,N,M2}. 
Again they also give an uniformization of the singular fibers of the family.
In particular the most degenerate fibers are uniformized by 
logarithmic derivatives of Schur polynomials. Those special properties of the
abelian functions are the reflection of the algebraic properties 
mentioned above. In the case of (n,s)-curves similar results are expected. 
Notice that this kind of degeneration structure from theta functions to 
trigonometric and rational functions are typical in the study of 
integrable systems.

Sato's theory of KP-hierarchy associates a point $U$ of the universal Grassmann
manifold (UGM) a so called tau function $\tau(t,\bxi)$, where $\bxi$ is a frame of $U$ \cite{S,SN}. It is a solution of the KP-hierarchy in the bilinear form. 
Conversely any solution of the KP-hierarchy can be written as a tau function of some point of $UGM$. It is
known that, for a point in a finite dimensional orbit of the KP-hierarchy,
the tau function can be expressed by the Riemann's theta function 
\cite{Mul,Kr,Shi,KNTY}. Moreover, 
for a frame $\bxi$, the expansion of the tau function can be explicitly
written as
\bea
&&
\tau(t,\bxi)=\sum \xi_\lambda s_\lambda(t),
\non
\ena
where the sum runs over the set of partitions $\lambda$,
$\xi_\lambda$ is the Pl\"ucker coordinate of $\bxi$ 
and $s_\lambda(t)$ is the Schur function corresponding to $\lambda$.
Thus, if the frame $\bxi$ is known, the expansion of 
the corresponding tau function is given very explicitly.

Now we consider the affine ring $A$ of an (n,s)-curve. 
A basis of $A$ as a vector space is given
explicitly by a set of certain monomials of $x$ and $y$. 
The space $A$ can be embedded in UGM using the local coordinate at $\infty$. 
We show that, for the normalized frame $\bxi^A$ of the corresponding point 
of UGM, the tau function is given by
\bea
&&
\tau(t;\bxi^A)=
\exp\left(-\sum_{i=1}^\infty c_it_i+\frac{1}{2}\qh(t)\right)\sigma(Bt),
\label{tau-sigma-intro}
\ena
where $\qh(t)=\sum \qh_{ij}t_it_j$, $B=(b_{ij})_{g\times \infty}$, $t={}^t(t_1,t_2,...)$ (see Theorems \ref{main} and \ref{tau-sigma}).
All constants $c_i$, $\qh_{ij}$, $b_{ij}$ are homogeneous 
polynomials of $\{\lambda_{kl}\}$ with rational coefficients. 
Using this formula we can deduce the above mentioned  
properties on the series expansion of a sigma function from 
those on the expansion of a tau function.

The crucial point to prove (\ref{tau-sigma-intro}) is the existence
of a holomorphic one form which vanishes at $\infty$ of order $2g-2$.
The square root of it plays the role of transforming half forms 
to functions \cite{N}. The existence of such a form is specific to (n,s)-curves.

Finally we remark that relations of tau and sigma functions are 
also discussed by C. Eilbeck, V. Enolski and J. Gibbons \cite{E,EEG}
extending the results of \cite{EH}.
In \cite{EEG} a similar relation to (\ref{tau-sigma-intro}) is derived.
The main difference between their formula and ours is that 
the frame $\bxi^A$ is described in terms of the derivatives 
of a tau function in the former
while it is described by the expansion coefficients of monomials
of $x$ and $y$ in the latter. They use their relation mainly to derive 
identities satisfied by Abelian functions.
While we use (\ref{tau-sigma-intro}) to give
a general algebraic formula for the expansion coefficients
of the tau function, the alternative to the sigma function. 
Therefore both results compensate to
each other and combining them is effective 
for a further development. 
Some related subjects are also studied in \cite{BS}.

The present paper is organized as follows.
After the introduction fundamental properties of an (n,s)-curve is explained
 in section 2.
In section 3 the KP-hierarchy and its reductions 
are reviewed. The Sato's theory of KP hierarchy and UGM is reviewed in section 4. In section 5 the embedding of the affine ring of an (n,s)-curve to UGM is
described. The construction of the sigma function is reviewed in section 7.
The expression of the tau function corresponding to the affine ring 
embedded in UGM is described in section 8. In section 9 the properties 
of the series expansion of the sigma function are studied based on the
formula in the previous section.

\section{$(n,s)$-Curve}
An $(n,s)$ curve is the plane algebraic curve defined by the equation 
$f(x,y)=0$ with
\bea
&&
f(x,y)=y^n-x^s-\sum_{ni+sj<ns}\lambda_{ij}x^iy^j,
\non
\ena
where $n,s$ are nonnegative integers which are coprime and satisfy 
$1<n< s$ \cite{BEL2,BL}. 
We assume it non-singular and denote by $X$ the corresponding
compact Riemann surface. Its genus is $g=1/2(n-1)(s-1)$.
Let $\pi:X\lar {\mathbb P}^1$ be the projection
to the $x$-coordinate, $(x,y)\lar x$. Then $\pi^{-1}(\infty)$ consists of
one point, which we denote by $\infty$, and is a branch point with the 
branching index $n$. 

The affine ring $A$ of $X$ is by definition
\bea
&&
A={\mathbb C}[x,y]/{\mathbb C}[x,y]f.
\non
\ena
Analytically it is isomorphic to the ring of meromorphic functions on $X$
which are holomorphic on $X-\{\infty\}$. For a meromorphic function $F$ on $X$
we denote $\ord\, F$ the order of poles at $\infty$. Then
\bea
&&
\ord\,x=n,
\qquad
\ord\,y=s.
\non
\ena

Let $\{f_i\}_{i=1}^\infty$ be the basis of $A$ as a vector space
specified by the conditions:
\vskip2mm
(i) $f_i\in \{x^{m_1}y^{m_2}\,|\,m_1\geq 0, n>m_2\geq 0\,\}$ ,

(ii) $\ord\, f_i<\ord\, f_{i+1}$ for $i\geq 1$.
\vskip5mm

\noindent
{\bf Example} $(n,s)=(2,2g+1)$: In this case $X$ is a hyperelliptic curve
of genus $g$. We have

$(f_1,f_2,...)=(1,x,x^2,...,x^{g},y,x^{g+1},xy,...)$.
\vskip5mm

Let $w_1<\cdots<w_g$ be the gap sequence at $\infty$. It is, by definition,
given by
\bea
&&
{\mathbb Z}_{\geq 0}\backslash \{\ord\, f_i\,|\,i\geq 1\,\}.
\non
\ena
In particular $w_1=1$ and $w_{g}=2g-1$.

\section{KP-hierarchy}

The following system of equations for a function $\tau(t)$, $t=(t_1,t_2,...)$
is called the bilinear equations of the KP-hierarchy \cite{DJKM}:
\bea
&&
\int_{k=\infty}\tau(t-[k^{-1}])\tau(t'+[k^{-1}])e^{\xi(t-t',k)}dk=0,
\non
\\
&&
\xi(t,k)=\sum_{i=1}^\infty t_i k^i,
\qquad
[k]=(k,\frac{k^2}{2},\frac{k^3}{3},...).
\label{bilinear-eq}
\ena

Notice that the equation (\ref{bilinear-eq}) is invariant under the 
the multiplication of the function of the form 
$c_0\exp(\sum_{i=1}^\infty c_i t_i)$ to $\tau(t)$, where $c_i$ are constants.

The equations (\ref{bilinear-eq}) can be rewritten using the Hirota derivative:
\bea
&&
D^M \tau(t)\cdot \tau(t)=\partial_y^M \tau(t+y)\tau(t-y)|_{y=0},
\non
\ena
where $D=(D_1,D_2,...)$, $y=(y_1,y_2,...)$, $M=(m_1,...,m_l)$, $l\geq 0$,
$D^M=D_1^{m_1}\cdots D_l^{m_l}$, $\partial_y^M=
\partial_{y_1}^{m_1}\cdots \partial_{y_l}^{m_l}$.
Let 
\bea
&&
e^{\xi(t,k)}=\sum_{j=0}^\infty p_j(t) k^j.
\non
\ena
Then (\ref{bilinear-eq}) is equivalent to
\bea
&&
\sum_{j=0}^\infty p_j(-2y)p_{j+1}(\tilde{D})e^{\sum_{l=1}^\infty y_l D_l}
\tau(t)\cdot \tau(t)=0,
\label{KP-Hirota}
\ena
where $\tilde{D}=(D_1,D_2/2,D_3/3,...)$ \cite{DJKM}.
The coefficient of $y_3$ gives the KP-equation in the bilinear form:
\bea
&&
(D_1^4+3D_2^2-4D_1D_3)\tau(t)\cdot \tau(t)=0.
\non
\ena

The system of equations obtained from (\ref{KP-Hirota}) by setting $D_{jn}=0$
for all $j\geq 1$ is called the n-reduced KP-hierarchy.
For example the bilinear form of the KdV equation
\bea
&&
(D_1^4-4D_1D_3)\tau(t)\cdot \tau(t)=0
\non
\ena
is the first member of the 2-reduced KP-hierarchy and
\bea
&&
(D_1^4+3D_2^2)\tau(t)\cdot \tau(t)=0
\non
\ena
is the first member of the 3-reduced KP-hierarchy (Boussinesq equation).

A solution $\tau(t)$ of the KP-hierarchy is a solution of the n-reduced
KP-hierarchy if $\tilde{\tau}(t)= c_0e^{\sum_{i=1}^\infty c_i t_i}\tau(t)$ 
does not depend on $\{t_{nj}\,|\,j\geq 1\}$ for some constants $c_i$.

\section{Universal Grassmann Manifold}
In this section we briefly review Sato's theory of KP equation
and the universal Grassmann manifold (UGM) following \cite{S,SN}
(see \cite{NT} for the English translation of \cite{SN}).

Let $R={\mathbb C}[[x]]$ be the ring of formal power series in $x$ and
${\cal E}_R=R((\partial^{-1}))$ the ring of microdifferential operators 
with the coefficients in $R$:
\bea
&&
{\cal E}_R=\{\sum_{-\infty<i<<\infty} a_i(x)\partial^i\,|\, a_i(x)\in R\},
\qquad
\partial=\frac{d}{dx}.
\non
\ena
Using the Leibnitz rule
\bea
&&
a(x)\partial^i=\sum_j (-1)^j\bc{i}{j}\partial^{i-j}a^{(j)}(x),
\qquad
a^{(j)}(x)=\frac{d^j a(x)}{d x^j},
\non
\ena
${\cal E}_R$ can be described as the set of 
operators of the form
$\sum_{-\infty<i<<\infty} \partial^i a_i(x)$ as well.

Let $V={\cal E}_R/{\cal E}_Rx\simeq {\mathbb C}((\partial^{-1}))$ be the left
${\cal E}_R$ module. We define the element $e_i$ of $V$ by 
\bea
&&
e_i=\partial^{-i-1}\quad \hbox{mod. ${\cal E}_Rx$}.
\non
\ena
The action of ${\cal E}_R$ on $e_i$ is given by
\bea
&&
\partial e_{i}=e_{i-1},
\qquad
x e_i=(i+1)e_i.
\non
\ena
We define two subspaces of $V$:
\bea
&&
V^\phi=\oplus_{i<0}{\mathbb C} e_i,
\qquad
V^{(0)}=\prod_{i>0}{\mathbb C} e_i.
\non
\ena
Then we have the decomposition
\bea
&&
V=V^\phi\oplus V^{(0)}.
\non
\ena
For a subspace $U$ of $V$ let 
\bea
&&
\pi_U: U\lar V/V^{(0)}\simeq V^\phi,
\non
\ena
be the composition of the inclusion $U \hookrightarrow V$ and the 
natural projection $V\lar V/V^{(0)}$.

\begin{defn}
The universal Grassmann manifold is the set of subspaces
$U$ of $V$ such that $\ker\, \pi_U$, $\coker\,\pi_U$ are
finite dimensional and the index of $\pi_U$ is zero:
\bea
&&
\ind(\pi_U)=\dim(\ker\, \pi_U)-\dim(\coker\,\pi_U)=0.
\non
\ena
\end{defn}

For a partition $\lambda=(\lambda_1,...,\lambda_l)$ we define
the Schur function $s_\lambda(t)$ by
\bea
&&
s_\lambda(t)=\det(p_{\lambda_i-i+j}(t))_{1\leq i,j\leq l}.
\non
\ena

For a point $U$ of UGM, a frame $\bxi$ of $U$ is a basis of $U$
\bea
&&
\bxi=(\xi_j)_{j<0}=(...,\xi_{-2},\xi_{-1}),
\quad
\xi_j=\sum_{i\in {\mathbb Z}} \xi_{ij} e_i,
\non
\ena
such that, for $j<<0$,
\bea
&&
\xi_{ij}=\left\{
\begin{array}{ll}
0&\quad i<j\\
1&\quad i=j.
\end{array}
\right.
\ena

For a point $U$ of UGM there exists a unique sequence of integers
$\rho_U=(\rho(i))_{i<0}$ and a unique frame $\bxi$ of $U$
such that
\bea
&&
\rho(-1)>\rho(-2)>\rho(-3)>\cdots,\qquad \rho(i)=i \hbox{ for $i<<0$},
\label{maya}
\\
&&
\xi_{ij}=\left\{
\begin{array}{ll}
0&\quad i<\rho(j) \mbox{ or $i=\rho(j')$ for some $j'>j$}\\
1&\quad i=\rho(j).
\end{array}
\right.
\label{n-frame}
\ena
The frame satisfying the condition (\ref{n-frame}) is said to be normalized.
 
In terms of $\rho_U=(\rho(i))_{i<0}$ the index of $\pi_U$
is given by
\bea
&&
\sharp\{i\,|\,\rho(i)\geq 0\,\}-
\sharp({\mathbb Z}_{<0}\backslash \{i\,|\,\rho(i)<0\,\}),
\non
\ena
where ${\mathbb Z}_{<0}$ is the set of negative integers.

In general a sequence of integers $\rho=(\rho(i))_{i<0}$
which satisfies (\ref{maya})
determines a partition 
$\lambda_{\bf \rho}=(\lambda_1,\lambda_2,...)$ by
\bea
&&
\lambda_i=\rho(-i)+i.
\non
\ena
\vskip5mm

\noindent
{\bf Example} If $\rho=(2,0,-2,-4,-5,-6,...)$ then
$\lambda_\rho=\rho-(-1,-2,-3,-4,...)=(3,2,1,0,...)=(3,2,1)$.
\vskip5mm

\begin{defn}
For a partition $\lambda$ the set of the points 
$U$ in $UGM$ satisfying $\lambda_\rho=\lambda$
is denoted by $UGM^\lambda$.
\end{defn}

The UGM is the disjoint union of $UGM^\lambda$'s:
\bea
&&
UGM=\coprod_{\lambda} UGM^\lambda.
\non
\ena

Given a point $U$ of $UGM$ and 
a sequence ${\bf \rho}=(\rho(i))_{i<0}$ satisfying the condition
(\ref{maya}), 
the Pl\"ucker coordinate $\xi_{\lambda_{\bf \rho}}$ of the normalized frame
$\bxi$
is defined as the determinant of the ${\mathbb Z}_{<0}\times {\mathbb Z}_{<0}$
matrix:
\bea
&&
\xi_{\lambda_\rho}=\det(\xi_{\rho(i)j})_{i,j<0}.
\non
\ena

For two partitions 
$\lambda=(\lambda_1,\lambda_2,...)$ and 
$\mu=(\mu_1,\mu_2,...)$, 
we define $\lambda\geq \mu$ if $\lambda_i\geq \mu_i$ for any $i$.

For $U\in UGM^\lambda$ the Pl\"ucker coordinates satisfy
\bea
&&
\xi_\mu=\left\{
\begin{array}{ll}
0&\quad \hbox{unless $\mu\geq \lambda$}\\
1&\quad \mu=\lambda.
\end{array}
\right.
\non
\ena

\begin{defn}
For $U\in UGM$ let
\bea
&&
\tau(t;\bxi)=\sum_{\mu\geq \lambda}\xi_\mu s_\mu(t),
\label{exp-tau}
\ena
where $\bxi$ is the normalized frame of $U$.
The function $\tau(t;\bxi)$ and its constant multiple is called the 
 tau function of $U$.
\end{defn}

\begin{theorem}\cite{S,SN}\label{UGM-tau}
For $U$ in $UGM^\lambda$ the tau function of $U$ is a solution of the
KP-hierarchy (\ref{bilinear-eq}). Conversely any solution $\tau(t)\in {\mathbb C}[[t]]$
of (\ref{bilinear-eq}) there exists a point of $U$ of UGM such that $\tau(t)$ is the tau
function of $U$.
\end{theorem}

The theorem follows from the fact \cite{SN,S,NT} that, if we expand $\tau(t)$ as
\bea
&&
\tau(t)=\sum_\lambda \xi_\lambda s_\lambda(t)
\non
\ena
with some set of constants $\{\xi_\lambda\}$, the bilinear equation 
(\ref{bilinear-eq}) is equivalent to the Pl\"ucker relations for 
$\{\xi_\lambda\}$. Based on this theorem, the point of UGM corresponding 
to the solution $\tau(t)$ of (\ref{bilinear-eq}) is recovered 
through the wave function as follows.

Let $K={\mathbb C}((x))$ be the field of formal Laurent series in $x$ and 
${\cal E}_K=K((\partial^{-1}))$ the ring of microdifferential operators:
\bea
&&
{\cal E}_K=\{\sum_{-\infty<i<<\infty}a_i(x)\partial^i\,|\,a_i(x)\in K\,\}.
\non
\ena

\begin{defn}\label{pseudo-reg}
Let ${\cal W}$ be the set of $W$ in ${\cal E}_K$ of the form
\bea
&&
W=\sum_{i\leq 0}w_i \partial^i,
\quad
w_0=1,
\non
\ena
satisfying the condition that there exist non-negative integers $l,m$ such that
\bea
&&
x^lW, \quad W^{-1}x^m\in {\cal E}_R.
\non
\ena
\end{defn}

Then 

\begin{theorem}\cite{S,SN}\label{W-UGM}
There is a bijective map $\gamma:{\cal W}\lar UGM$ given by
\bea
&&
\gamma(W)=W^{-1}x^mV^\phi,
\non
\ena
where $m$ is chosen as in Definition \ref{pseudo-reg}. 
The image $\gamma(W)$ does not depend on the choice of $m$.
\end{theorem}

Let $\tau(t)\in {\mathbb C}[[t]]$ be a solution of the KP-hierarchy 
(\ref{bilinear-eq}).
The wave function and the conjugate wave function are defined by
\bea
&&
\Psi(t;z)=\frac{\tau(t+[z])}{\tau(t)}\exp(-\sum_{i=1}^\infty t_iz^{-i}),
\non
\\
&&
\bar{\Psi}(t;z)=\frac{\tau(t-[z])}{\tau(t)}\exp(\sum_{i=1}^\infty t_iz^{-i}).
\non
\ena
Due to the bilinear identity (\ref{bilinear-eq}) there exists $W\in {\cal W}$
such that $\Psi$ and $\bar{\Psi}$ can be written as \cite{DJKM}
\bea
&&
\Psi(t,z)=(W^{\ast})^{-1}\exp(-\sum_{i=1}^\infty t_iz^{-i}),
\non
\\
&&
\bar{\Psi}(t,z)=W\exp(\sum_{i=1}^\infty t_iz^{-i}),
\non
\ena
where
$x=t_1$ and $P^\ast=\sum (-\partial)^i a_i(x)$ is 
the formal adjoint of $P=\sum a_i(x) \partial^i$.
We have $(P^\ast)^{-1}=(P^{-1})^\ast$ 
for an invertible $P\in {\cal E}_K$.

The following lemma easily follows from the definition of $s_\lambda(t)$.

\begin{lemma}\cite{SW} 
For any partition $\lambda=(\lambda_1,\lambda_2,...)$ 
\bea
&&
s_\lambda(t_1,0,0,...)=d_\lambda t_1^{|\lambda|},
\non
\\
&&
d_\lambda=\frac{\prod_{i<j}(\mu_i-\mu_j)}{\prod_{i=1}^l\mu_i!},
\non
\ena
where $l$ is taken large enough such that $\lambda_l=0$,
$\mu_i=\lambda_i+l-i$ and $|\lambda|=\sum_i \lambda_i$.
\end{lemma}

By the lemma and Theorem \ref{UGM-tau} we see that $\tau(x,0,0,...)$ 
is not identically zero. 
Let $m_0$ be the order of zeros of $\tau(x,0,0...)$ at $x=0$ 
and $m\geq m_0$. Obviously we have 
$x^m W(x,0,...), W(x,0,...)^{-1}x^m\in {\cal E}_R$ which implies
$W\in {\cal W}$.
Then 

\begin{theorem}\cite{SN,S}
\bea
&&
\tau(t)=\tau\left(t;\gamma\left(W(x,0,...)\right)\right).
\non
\ena
\end{theorem}

Let us describe $\gamma(W(x,0,...))$ in terms of the wave function $\Psi$.
\begin{prop}
Let
\bea
&&
x^m\Psi(x,0,...;z)=\sum_{i=0}^\infty \Psi_i(z)\frac{x^i}{i!}.
\non
\ena
Then we have
\bea
&&
(-1)^iW(x,0,...)^{-1}x^me_{-1-i}=\Psi_i(\partial^{-1})e_{-1},
\non
\\
&&
\gamma(W(x,0,...))=\mbox{Span}_{\mathbb C}\{\Psi_i(\partial^{-1})e_{-1}\,|\,
i\geq 0\,\},
\non
\ena
where $\mbox{Span}_{\mathbb C}\{\cdots\}$ signifies the vector space generated
by $\{\cdots\}$.
\end{prop}
\vskip2mm
\noindent
{\it Proof.}
Let 
\bea
&&
x^m\Psi(x,0,...;z)=\sum_{i=0}^\infty \frac{x^i}{i!}\psi_i(z)e^{-xz^{-1}}
=\sum_{i=0}^\infty \frac{x^i}{i!}\psi_i(-\partial^{-1})e^{-xz^{-1}}.
\non
\ena
Then
\bea
&&
W(x,0...)^{-1}x^m=
\left(\sum_{i=0}^\infty \frac{x^i}{i!}\psi_i(-\partial^{-1})\right)^{\ast}
=\sum_{i=0}^\infty \psi_i(\partial^{-1})\frac{x^i}{i!}.
\non
\ena
Thus
\bea
W(x,0,...)^{-1}x^me_{-1-i}
&=&
\sum_{j=0}^\infty \psi_j(\partial^{-1})\frac{x^j}{j!}
\partial^ie_{-1}
\non
\\
&=&
\sum_j (-1)^j\bc{i}{j}\psi_j(\partial^{-1})\partial^{i-j}e_{-1}.
\label{cal-1}
\ena
On the other hand
\bea
\Psi_i(z)&=&
\frac{\partial^i}{\partial x^i}(x^m\Psi(x,0,...|z))\vert_{x=0}
\non
\\
&=&
\sum_{j=0}^i\bc{i}{j}\psi_j(z)\left(-\frac{1}{z}\right)^{i-j}.
\label{cal-2}
\ena
The assertion of the lemma follows from (\ref{cal-1}), (\ref{cal-2})
and Theorem \ref{W-UGM}.
\qed

\section{Embedding of the Affine Ring to UGM}

We can take a local coordinate $z$ of $X$ around $\infty$ in such a way
that
\bea
x=\frac{1}{z^n},
\qquad
y=\frac{1}{z^s}(1+O(z)).
\label{coord}
\ena

Using the expansion in $z$ define the embedding $\iota:A\lar V$ by
\bea
&&
\sum a_m z^m \mapsto \sum a_m e_{m+g-1}.
\non
\ena

Let $U^A=\iota(A)$.
\begin{lemma}
The image $U^A$ belongs to $UGM$.
\end{lemma}
\vskip2mm
\noindent
{\it Proof.} Let $0=w_1^\ast<w_g^{\ast}<\cdots$ be non-gaps 
of $X$ at $\infty$, that is, $\{w_i^\ast\}=\{\ord\,f_i\}$. Then
They satisfy
\bea
&&
w_{i}^\ast=g-1+i\quad\hbox{for $i\geq g+1$},
\non
\\
&&
\{\,w_1^\ast,...,w_g^\ast\,\}\sqcup\{w_1,...,w_g\}=\{0,1,...,2g-1\}.
\label{gap-nongap}
\ena
We have
\bea
&&
\dim\,(\ker\,\pi_{U^A})=\sharp\{\,i\,|\,w_i^\ast\leq g-1\},
\quad
\dim\,(\coker\,\pi_{U^A})=\sharp\{\,i\,|\,w_i> g-1\}.
\non
\ena
Then the equation (\ref{gap-nongap}) implies that
\bea
&&
(g-\dim\,(\ker\,\pi_{U^A}))+\dim\,(\coker\,\pi_{U^A})=g.
\non
\ena
Thus the index of $\pi_{U^A}$ is zero and $U^A$ is in UGM.
\qed

Let
\bea
&&
\lambda(n,s)=(w_g,...,w_1)-(g-1,...,1,0)
\non
\ena
be the partition associated with the gap sequence.
Then $U^A$ belongs to $UGM^{\lambda(n,s)}$ and
the tau function $\tau(t;\bxi^A)$ has the expansion
\bea
&&
\tau(t;\bxi^A)=s_{\lambda(n,s)}(t)+
\sum_{\lambda>\lambda(n,s)}\xi^A_\lambda s_\lambda(t),
\label{exp-tau-ua}
\ena
where $\bxi^A$ is the normalized frame of $U^A$.
The aim of the paper is to determine the analytic expression 
of $\tau(t;\bxi^A)$.

\section{Sigma Functions}

Let $X$ be an $(n,s)$ curve introduced in section 2
and $\{du_{w_i}\}$ a basis of holomorphic one forms given by
\bea
&&
du_{w_i}=-\frac{f_{g+1-i}dx}{f_y},
\non
\ena where $f_i$ is the monomial of $x,y$ defined in section 2.
We choose an algebraic fundamental form $\homega(p_1,p_2)$
\cite{N} and decompose it as
\bea
&&
\homega(p_1,p_2)=d_{p_2}\Omega(p_1,p_2)+\sum_{i=1}^g du_{w_i}(p_1)dr_i(p_2),
\non
\ena
where $dr_i(p)$ is a second kind differential holomorphic outside $\infty$.
Here
\bea
&&
\Omega(p_1,p_2)=
\frac{\sum_{i=0}^{n-1}y_1^i[\frac{f(z,w)}{w^{i+1}}]_{+}\vert_{(z,w)=(x_2,y_2)}}
{(x_1-x_2)f_y(x_1,y_1)}dx_1,
\non
\\
&&
[\sum_{m\in {\mathbb Z}} a_m w^m]_{+}=\sum_{m\geq 0} a_m w^m.
\non
\ena
Then, with respect to the intersection form $\circ$ defined by
\bea
&&
\omega\circ \eta=\Res_{p=\infty} (\int^p\omega)\eta.
\non
\ena
$\{du_{w_i},dr_j\}$ becomes symplectic:
\bea
&&
du_{w_i}\circ du_{w_j}=dr_i\circ dr_j=0, 
\qquad
du_{w_i}\circ dr_j=\delta_{ij}.
\label{symp}
\ena

Let us take a symplectic basis $\{\alpha_i,\beta_j\}$ of the homology group 
of $X$ and form the period matrices
\bea
&&
2\omega_1=(\int_{\alpha_j}du_{w_i}),
\qquad
2\omega_2=(\int_{\beta_j}du_{w_i}),
\non
\\
&&
-2\eta_1=(\int_{\alpha_j}dr_i),
\qquad
-2\eta_2=(\int_{\beta_j}dr_i),
\non
\\
&&
\tau=\omega_1^{-1}\omega_2.
\ena


Let $\delta=\delta'\tau +\delta''$, $\delta',\delta''\in {\mathbb R}^g$ 
be a representative of Riemann's constant with respect to the choice
$(\{\alpha_i,\beta_j\},\infty)$ and set $\delta={}^t(\delta',\delta")
\in {\mathbb R}^{2g}$.

In general, for $a,b\in {\mathbb R}^g$  and a point $\tau$ of 
the Siegel upper half space of degree $g$, 
the Riemann's theta function is defined by
\bea
&&
\theta\qbc{a}{b}(z,\tau)=\sum_{m\in {\mathbb Z}^g}
\exp(\pi i{}^t(m+a)\tau (m+b)+2\pi i{}^t (m+a)(z+b)),
\non
\ena
where $z={}^t(z_1,...,z_g)$ \cite{M1}.

\begin{defn}
We define
\bea
&&
{\hat \sigma}(u)
=\exp(\frac{1}{2}{}^tu\eta_1\omega_1^{-1}u)
\theta[\delta]((2\omega_1)^{-1}u,\tau),
\non
\ena
where $u={}^t(u_{w_1},...,u_{w_g})$.
\end{defn}

Notice that ${\hat \sigma}(u)$ depends on the choice of 
$\{du_{w_i}\}$, $\homega(p_1,p_2)$,
and $\{\alpha_i,\beta_j\}$. For the degrees of freedom on the choice 
$\homega(p_1,p_2)$ see section 3.4 of \cite{N}.

Later we define the sigma function $\sigma(u)$ by multiplying
a certain constant to ${\hat \sigma}(u)$ (Definition \ref{sigma}).

The function $\hat{\sigma}(u)$ has the following transformation rule 
(see \cite{N}).

\begin{prop}\cite{BEL1}
For $m_1,m_2\in {\mathbb Z}^g$ 
\bea
&&
\hat{\sigma}(u+2\omega_1m_1+2\omega_2m_2)
\non
\\
&&
=(-1)^{{}^tm_1m_2+2{}t(\delta'm_1-\delta''m_2)}
\exp\left({}^t(2\eta_1m_1+2\eta_2m_2)(u+\omega_1m_1+\omega_2m_2)\right)
\hat{\sigma}(u).
\non
\ena
\end{prop}

\section{Tau Function in Terms of Sigma Function}
Let us take the local coordinate $z$ of $X$ around $\infty$ as 
in (\ref{coord}) and
consider the expansions in $z$:
\bea
du_{w_i}&=&\sum_{j=1}^\infty b_{ij}z^{j-1}dz,
\non
\\
\homega(p_1,p_2)&=&
\left(
\frac{1}{(z_1-z_2)^2}+\sum_{i,j\geq 1}\qh_{ij}z_1^{i-1}z_2^{j-1}
\right)
dz_1dz_2,
\non
\ena
By the definition of $w_i$, $du_i$ and $z$ we have \cite{N}
\bea
&&
b_{ij}=\left\{
\begin{array}{cc}
0& \quad\hbox{if $j<w_i$}\\
1& \quad\hbox{if $j=w_i$}.
\end{array}
\right.
\label{b-ij}
\ena

In particular we have
\bea
&&
du_{w_g}=z^{2g-2}(1+\sum_{j>2g-2}b_{gj}z^{j-2g+1})dz.
\non
\ena
Let 
\bea
&&
\log z^{-(g-1)}\sqrt{\frac{du_{w_g}}{dz}}=
\sum_{i=1}^\infty \frac{c_i}{i}z^i,
\non
\\
&&
B=(b_{ij})_{g\times \infty},
\qquad
t={}^t(t_1,t_2,...),
\non
\\
&&
\qh(t)=\sum_{i,j=1}^\infty \qh_{ij}t_it_j.
\non
\ena
Then

\begin{theorem}\label{main}
(i) There exists a constant $C$ such that
\bea
&&
\tau(t;\bxi^A)=
C \exp\left(-\sum_{i=1}^\infty c_it_i+\frac{1}{2}\qh(t)\right)\hat{\sigma}(Bt).
\label{tau=sigma}
\ena
\vskip2mm
\noindent
(ii) The tau function $\tau(t;\bxi^A)$ satisfies the $n$-reduced KP-hierarchy.
\end{theorem}
\vskip2mm
\noindent
{\t Proof.} We denote the right hand side of (\ref{tau=sigma}) by 
$\widehat{\tau}(t)$

Let $E(p_1,p_2)$ be the prime form \cite{F}. 
Using the local coordinate $z$ we define $E(z_1,z_2)$, $E(\infty,p)$ by
\bea
&&
E(p_1,p_2)=\frac{E(z_1,z_2)}{\sqrt{dz_1}\sqrt{dz_2}},
\non
\\
&&
E(\infty,p_2)=\frac{E(0,z_2)}{\sqrt{dz_2}}.
\non
\ena
The normalized fundamental form $\omega(p_1,p_2)$
is defined by
\bea
&&
\omega(p_1,p_2)=d_{p_1}d_{p_2}\log\,E(p_1,p_2).
\non
\ena
It has the expansion of the form
\bea
&&
\omega(p_1,p_2)=\left(\frac{1}{(z_1-z_2)^2}+\sum_{i,j\geq 1}q_{ij}z_1^{i-1}z_2^{j-1}\right)dz_1dz_2.
\non
\ena
Let $\{dv_j\}_{j=1}^g$ be the normalized basis of holomorphic 1-forms and
\bea
&&
dv_i=\sum_{j=1}^\infty a_{ij} z^{j-1}
\non
\ena
its expansion near $\infty$. We set
\bea
&&
\bar{A}=(a_{ij}),
\qquad
q(t)=\sum_{i,j=1}^\infty q_{ij}t_it_j.
\non
\ena
Then the following theorem is well known.

\begin{theorem}\cite{Kr,Shi}
The function
\bea
&&
\tilde{\tau}(t)=\exp\left(\frac{1}{2}q(t)\right)\theta(\bar{A}t+\zeta)
\non
\ena
is a solution of the KP-hierarchy for any $\zeta\in {\mathbb C}^g$.
\end{theorem}

Taking $\zeta=\delta$, using the relations \cite{N}
\bea
&&
2\omega_1\bar{A}=B,
\non
\\
&&
\homega(p_1,p_2)=\omega(p_1,p_2)-
\sum_{i,j=1}^g(\eta_1\omega_1^{-1})_{ij}du_{w_i}(p_1)du_{w_j}(p_2),
\label{omega-homega}
\ena
and the definition of $\sigma$ in terms of the Riemann's theta we easily
see that $\widehat{\tau}(t)$ is obtained from $\tilde{\tau}(t)$ by
multiplying a constant and the exponential of a linear function of $t$. Thus
$\widehat{\tau}(t)$ is a solution of the KP-hierarchy.

In order to determine the point of UGM corresponding to $\widehat{\tau}(t)$
we calculate the wave function. 

Let $d\tilde{r}_i$ be the normalized abelian differential of the second
kind which means that it is holomorphic on $X-\{\infty\}$, has zero $\alpha_j$
period for any $j$ and it has the form 
\bea
&&
d\tilde{r}_i=d\left(\frac{1}{z^i}+O(1)\right),
\non
\ena
near $\infty$. We set
\bea
&&
d\widehat{r}_i=
d\tilde{r}_i+\sum_{k,l=1}^g b_{ki}(\eta_1\omega_1^{-1})_{kl}du_{w_l}.
\non
\ena

By calculation we have
\bea
&&
\Psi(t,z)
\non
\\
&=&
\frac{\widehat{\tau}(t+[z])}{\widehat{\tau}(t)}\exp(-\sum_{i=1}^\infty t_iz^{-i})
\non
\\
&=&\sqrt{\frac{du_{w_g}}{dz}}^{-1}\frac{z^g}{E(0,z)}
\frac{\sigma(Bt+\int_\infty^p du)}{\sigma(Bt)}
\exp\left(-\sum_{i=1}^\infty t_i\int^p d\widehat{r}_i
-\frac{1}{2}\int_\infty^p {}^tdu \cdot \eta_1(\omega_1)^{-1} \cdot 
\int_\infty^p du
\right).
\non
\ena
This is simply a restatement of the known result (\cite{Kr}, \cite{KNTY}).

The following modification of the prime form had been introduced in \cite{N}:
\bea
&&
\tilde{E}(\infty,p)=E(\infty,p)\sqrt{du_{w_g}}
\exp(\frac{1}{2}\int_\infty^p {}^tdu\cdot \eta_1(\omega_1)^{-1}\cdot
 \int_\infty^p du).
\non
\ena
Notice that this is not a half form but a multi-valued holomorphic function 
on $X$ with zeros only at $\infty$. The expansion near $\infty$ is of the form
\bea
&&
\tilde{E}(\infty,p)=z^g(1+O(z)).
\non
\ena
Its transformation rule is determined in \cite{N}.

\begin{lemma}\cite{N}
Let $\gamma\in \pi_1(X,\infty)$. Suppose that its image in $H_1(X,{\mathbb Z})$
is given by $\sum_{i=1}^gm_{1,i}\alpha_i+\sum_{i=1}^gm_{2,i}\beta_i$.
Then 
\bea
&&
\tilde{E}(\infty,\gamma(p))/\tilde{E}(\infty,p)
\non
\\
&&
=(-1)^{{}^tm_1m_2+2({}^t\delta'm_1-{}^t\delta''m_2)}
\exp\left(
{}^t(2\eta_1m_1+2\eta_2m_2)(\int_\infty^p du+\omega_1m_1+\omega_2m_2)
\right),
\non
\ena
where $m_i={}^t(m_{i,1},...,m_{i,g})$.
\end{lemma}

We rewrite $\Psi$ using $\tilde{E}(\infty,p)$ as
\bea
&&
\Psi(t,z)=
\frac{z^g}{\tilde{E}(\infty,p)}
\frac{\sigma(Bt+\int_\infty^p du)}{\sigma(Bt)}
\exp\left(-\sum_{i=1}^\infty t_i\int^p d\widehat{r}_i
\right).
\non
\ena

\begin{lemma}
The function $z^{-g}\Psi(t,z)$ is $\pi_1(X,\infty)$-invariant 
and any coefficient of $t_1^{m_1}t_2^{m_2}\cdots$ in
the expansion of $\sigma(Bt)z^{-g}\Psi(t,z)$ is in $A$.
\end{lemma}

This lemma follows from 

\begin{lemma}
We have
\bea
&&
\int_{\alpha_j}d\widehat{r}_i=\left({}^t(2\eta_1)B\right)_{ij},
\qquad
\int_{\beta_j}d\widehat{r}_i=\left({}^t(2\eta_2)B\right)_{ij}.
\non
\ena
\end{lemma}

This lemma can be proved by a direct calculation using the definition of
$d\tilde{r}_i$.

Since $\widehat{\tau}(t)$ is a tau function of the KP-hierarchy,
$\widehat{\tau}(x,0,...)$ is not identically zero. 
Let $m$ be the order of zeros of $\widehat{\tau}(x,0,...)$ at $x=0$ and
\bea
&&
x^m\Psi(x,0,...;z)=\sum_{i=0}^\infty \Psi_i(z)x^i.
\non
\ena
Then 
\bea
&&
z^{-g}\Psi_i(z)\in A.
\non
\ena
Let
\bea
&&
z^{-g}\Psi_i(z)=\sum_{-\infty<<k<\infty}\psi_k z^k.
\non
\ena
Then
\bea
&&
\Psi_i(\partial^{-1})e_{-1}=\sum_k \psi_k e_{g-1+k}=\iota(z^{-g}\Psi_i(z)).
\non
\ena
Thus the subspace $U$ of $V$ generated by $\{\Psi(\partial^{-1})e_{-1})\}$
is a subspace of $U^A$. Since both $U$ and $U^A$ are in UGM, $U=U^A$.

Next we prove that $\tau(t;\bxi^A)$ is a solution of the n-reduced
KP-hierarchy.

\begin{lemma}\label{n-reduced}
\bea
&&
\qh_{ij}=0 \quad \hbox{if $i$ or $j=0$ mod.$n$},
\non
\\
&&
b_{ij}=0 \quad \hbox{if $j=0$ mod.$n$}.
\non
\ena
\end{lemma}
\vskip2mm
\noindent
{\it Proof.} Firstly let us prove $b_{ij}=0$ if $j=0$ mod.$n$.
Notice that
\bea
&&
d\tilde{r}_{nk}=dx^k=d\left(\frac{1}{z^{nk}}\right).
\label{exp-secondkind}
\ena
Then
\bea
&&
d\tilde{r}_{nk}\circ du_{w_i}
=\Res_{p=\infty}\left((\int^p d\tilde{r}_{nk})du_{w_i}\right)=b_{i,nk}.
\non
\ena
On the other hand the left hand side is zero because $d\tilde{r}_{nk}$ is an
exact form.

Next we prove $q_{ij}=0$ if $i,j$ satisfy $ij=0$ mod. $n$. 
In fact $q_{ij}$ can be obtained
as the expansion of $d\tilde{r}_{i}$ as (see the appendix of \cite{KNTY} for
example)
\bea
&&
d\tilde{r}_{i}=d(\frac{1}{z^i}-\sum_{j=1}^\infty q_{ij}\frac{z^j}{j}).
\non
\ena
Then the assertion follows from (\ref{exp-secondkind}) and the symmetry of
$q_{ij}$.

Finally $\qh_{ij}=0$ for $(i,j)$ satisfying $ij=0$ mod.$n$ follows 
from the relation (\ref{omega-homega}) and the symmetry of $\qh_{ij}$.
\qed

It follows from the lemma that $\exp(\sum_{i\geq 1}c_it_i)\tau(t;\bxi^A)$ 
does not depend on
$t_{nk}$ for $k\geq 1$. Thus $\tau(t;\xi^A)$ is a solution
of the n-reduced KP-hierarchy.
\qed

\begin{cor} The coefficients of $x^l$, $l\geq 0$, in the expansion 
of the function
\bea
&&
\frac{\sigma((x,0,...)+\int_\infty^p du)}{\tilde{E}(\infty,p)}
\exp\left(- x\int^p d\widehat{r}_1\right)
\non
\ena
generate the affine ring $A$ as a vector space.
\end{cor}

\vskip2mm
\noindent
{\bf Remark.} In the case of $g=1$ the above corollary tells that
the coefficients of $x^l$, $l\geq 0$ of the Baker-Akhiezer function 
\bea
&&
\frac{\sigma(u+x)}{\sigma(u)}\exp\left(-x\zeta(u)\right)
\non
\ena
generate the space generated by $\wp^{(i)}(u)$, $i\geq 0$ and $1$.
This fact can be easily checked and is well known.

\section{Applications}

In this section we study the series expansion of the sigma function
as an application of Theorem \ref{main}.

Let $u={}^t(u_{w_1},...,u_{w_g})$.
We define the degrees of $u_i$ and $t_i$ to be $-i$:
\bea
&&
\deg\,u_i=\deg\,t_i=-i.
\non
\ena

\begin{theorem}\label{normalization}
For the constant $C$ in Theorem \ref{main} we have the following 
series expansion at $u=0$:
\bea
&&
C \theta[\delta]((2\omega_1)^{-1}u,\tau)=s_{\lambda(n,s)}(u)+\cdots,
\non
\ena
where $\cdots$ part contains lower degree terms than $s_{\lambda(n,s)}(u)$.
\end{theorem}
\vskip2mm
\noindent
{\it Proof.} Let $t^0$ be $t=(t_1,t_2,...)$ in which $t_j=0$ for $j\notin 
\{w_1,...,w_g\}$ and $u=Bt^0$. Then, by (\ref{b-ij}), 
we have
\bea
&&
u_{w_i}=t_{w_i}+\sum_{j=i+1}^g b_{i w_j}t_{w_j}.
\non
\ena
Inverting this we have
\bea
&&
t_{w_i}=u_{w_i}+\sum_{j=i+1}^g b'_{i w_{j}}u_{w_j}.
\label{t-u}
\ena
Then the theorem follows from Theorem \ref{main} and the expansion 
(\ref{exp-tau-ua}) of $\tau(t;\xi^A)$.
\qed

\begin{defn}\label{sigma}
The sigma function associated to the choice 
$(\{du_{w_i}\},\homega,\{\alpha_i,\beta_j\})$ is defined
by
\bea
&&
\sigma(u)=\sigma(u|\{du_i\},\homega(p_1,p_2),\{\alpha_i,\beta_j\})
=C\hat{\sigma}(u),
\non
\ena
where $C$ is the constant given in Theorem \ref{normalization}.
\end{defn}

The sigma function inherits some remarkable properties from 
the tau function $\tau(t;\bxi^A)$,
the algebraic expansion and the modular invariance.

The symplectic group $Sp(2g,{\mathbb Z})$ acts on the set of
canonical homology bases by
\bea
&&
M\left(\begin{array}{c}\alpha\\\beta\end{array}\right)
=\left(\begin{array}{cc}D&C\\B&A\end{array}\right)
\left(\begin{array}{c}\alpha\\\beta\end{array}\right),
\quad
\hbox{for } 
M=\left(\begin{array}{cc}A&B\\C&D\end{array}\right)
\in \hbox{Sp}(2g,{\mathbb Z}),
\non
\ena
where 
\bea
&&
\left(\begin{array}{c}\alpha\\\beta\end{array}\right)=
{}^t(\alpha_1,..,\alpha_g,\beta_1,...,\beta_g).
\non
\ena

We assign degrees to the coefficients of $f(x,y)$ as
\bea
&&
\deg\,\lambda_{ij}=ns-ni-sj.
\non
\ena

Then

\begin{theorem}\label{app-2}
(i) At $u=0$ $\sigma(u)$ has the following expansion:
\bea
&&
\sigma(u)=s_{\lambda(n,s)}(u)+\sum a_{\gamma} u^\gamma,
\non
\ena
where $\gamma=(\gamma_1,...,\gamma_g)$, 
$u^\gamma=u_{w_1}^{\gamma_1}\cdots u_{w_g}^{\gamma_g}$, $a_\gamma$ is
a homogeneous polynomial of $\{\lambda_{ij}\}$ with the degree
$-|\lambda(n,s)|+\sum_{i=1}^g \gamma_i w_i$ and the summation is taken
for $\gamma$ with $\sum_{i=1}^g \gamma_i w_i>|\lambda(n,s)|$.
\vskip2mm
\noindent
(ii) For $M\in \hbox{Sp}(2g,{\mathbb Z})$
\bea
&&
\sigma(u|\{du_i\},\homega,M{}^t(\alpha,\beta))=
\sigma(u|\{du_i\},\homega,{}^t(\alpha,\beta)).
\non
\ena
\end{theorem}

Notice that the property (i) implies the property (ii).
It is possible to study the modular transformation of the sigma function
using that of Riemann's theta function. However it is difficult to determine 
the 8-the root of unity part in that calculation. 

We remark that Theorem \ref{normalization} and \ref{app-2} had been
proved in \cite{N} in a different way.

\vskip2mm
\noindent
{\it Proof of Theorem \ref{app-2}.} 
\vskip2mm
\noindent
By Lemma 15 of \cite{N} $b_{ij}$ is a homogeneous polynomial of 
$\{\lambda_{kl}\}$
of degree $-w_i+j$ with the coefficient in ${\mathbb Q}$. In particular 
$\deg\,b_{gj}=-(2g-1)+j$. It follows that $c_i$ belongs to 
${\mathbb Q}[\{\lambda_{kl}\}]$ and it is homogeneous of degree $-i$.
We assign degree $-1$ to $z$. Then $\sum_{i\geq 1} c_it_i$ is homogeneous
of degree $0$. Also $\widehat{q}_{ij}$ belongs to 
${\mathbb Q}[\{\lambda_{kl}\}]$ and it is homogeneous of degree $i+j$
by Lemma 15 of \cite{N}. Thus $\deg\,\qh(t)=0$.

Let us calculate the degree of the Pl\"ucker coordinate $\xi^A_\lambda$.
Let $\lambda$ and $\lambda(n,s)$ correspond to $(\rho(i))_{i<0}$
and $(\rho_0(i))_{i<0}$ respectively and $0=w_1^\ast<w_2^\ast<\cdots$
the non-gaps of $X$ at $\infty$. Then
\bea
&&
\rho_0(-i)=-w_i^\ast+g-1,
\quad
\lambda(n,s)_i=\rho_0(-i)+i.
\label{rho-0}
\ena
Let us expand $f_j$ as
\bea
&&
f_j=\sum_{-\infty<<j<\infty}a_{ij}z^i.
\non
\ena
Then
\bea
&&
\xi^A_j=\iota(f_j)=\sum a_{i-g+1,j}e_i.
\non
\ena
Thus
\bea
&&
\xi^A_{ij}=a_{i-g+1,j}.
\non
\ena
By Lemma 15 in \cite{N} $a_{ij}\in {\mathbb Q}[\{\lambda_{kl}\}]$ 
and it is homogeneous of degree $i+w_j^\ast$. Therefore 
\bea
&&
\deg\,\xi^A_{ij}=i-g+1+w_j^\ast.
\non
\ena
Recall that
\bea
\xi^A_\lambda&=&\det(\xi^A_{\rho(-i),-j})_{i,j\geq 1}
\non
\\
&=&
\sum_{\sigma\in S_m}\xi^A_{\rho(-\sigma_1),-1}\cdots \xi^A_{\rho(-\sigma_m),-m},\non
\ena
where $m$ is taken large enough so that $\rho(i)=i$ for $i<-m$.
We have
\bea
\deg\,\xi^A_{\rho(-\sigma_1),-1}\cdots \xi^A_{\rho(-\sigma_m),-m}
&=&
\sum_{i\geq 1}(\rho(-i)-g+1+w_i^\ast)
\non
\\
&=&
\sum_{i\geq 1}(\rho(-i)-\rho_0(-i))
\non
\\
&=&
|\lambda|-|\lambda(n,s)|.
\non
\ena
Here we use (\ref{rho-0}) to eliminate $w_i^\ast$. Consequently
 $\xi^A_\lambda$ is a homogeneous polynomial of $\{\lambda_{ij}\}$
with the coefficients in ${\mathbb Q}$ of degree $|\lambda|-|\lambda(n,s)|$.
In the expression (\ref{t-u}) $b_{ij}'$ has the same properties as $b_{ij}$.
Thus the theorem follows from Theorem \ref{main} and (\ref{exp-tau-ua}).
\qed

We rephrase Theorem \ref{main} (i) in terms of the sigma function:

\begin{theorem}\label{tau-sigma}
\bea
&&
\tau(t;\bxi^A)=
\exp\left(-\sum_{i=1}^\infty c_it_i+\frac{1}{2}\qh(t)\right)\sigma(Bt).
\non
\ena
\end{theorem}

\vskip10mm

\noindent
{\large {\bf Acknowledgements}} 
\vskip3mm
\noindent
I would like to thank Masatoshi Noumi and Takashi Takebe 
for permitting me to see the manuscript \cite{NT} prior to its publication. 
I am also grateful to Yasuhiko Yamada for valuable comments.
Finally I would like to thank Victor Enolski for the explanation of their works
and useful comments on this paper. 
This research is supported by Grant-in-Aid for Scientific Research (B) 17340048.

\end{document}